\documentclass[a4paper,12pt,twoside]{article}
\usepackage{amsthm,amssymb,amsmath,amscd}
\usepackage{pb-diagram}
\newtheorem{thm}{Theorem}
\newtheorem{conj}{Conjecture}

\newtheorem{lem}{Lemma}

\newtheorem{defn}{Definition}

\newcommand{\Zset}{\mathbb{Z}}

\newcommand{\Cset}{\mathbb{C}}

\newcommand{\bc}{\mathbf{c}}

\def\Cent{{\rm Z}}

\def\Sym{{\rm Sym}}

\def\rec{{\rm rec}}
\def\red{{\rm red}}
\def\ssimple{{\rm ss}}

\def\HP{{\rm HP}}
\def\SL{{\rm SL}}

\def\SU{{\rm SU}}
\def\GL{{\rm GL}}

\def\Hom{{\rm Hom}}

\def\pt{{\rm pt}}

\def\Nor{{\rm N}}

\def\Irr{{\rm Irr}}

\def\Prim{{\rm Prim}}

\def\fs{{\mathfrak s}}

\def\fB{{\mathfrak B}}

\def\fB{{\mathfrak B}}

\def\q{{/\!/}}
\begin{document}
\title{Geometric structure in the representation theory of $p$-adic groups}
\author{Anne-Marie Aubert, Paul Baum and Roger Plymen}
\date{}
\maketitle

\begin{abstract}
We conjecture the existence of a simple geometric structure
underlying questions of reducibility of parabolically induced
representations of reductive $p$-adic groups.\end{abstract}

\section{Introduction}
\label{intro} In the representation theory of reductive $p$-adic
groups, the issue of reducibility of induced representations is an
issue of great intricacy: see, for example, the classic article by
Bernstein-Zelevinsky \cite{BZ1} on $\GL(n)$ and the more recent
article by Mui\'c \cite{M} on $G_2$.
%involving the $R$-group, Plancherel
%measure, local Langlands conjecture,  $L$-packets, endoscopy, etc.
It is our contention, expressed as a conjecture, that there exists
a simple geometric structure underlying this intricate theory.

For the moment, our conjecture is \emph{local}, in that it applies
only to finite places.  To explain our conjecture, we need to
refine the usual concept of quotient space.

\section{The extended quotient and the reduced quotient}
\label{extended quotient} Let $\Gamma$ be a finite group and let
$X$ be a complex affine algebraic variety. Assume that $\Gamma$ is
acting on $X$ as automorphisms of the affine algebraic variety
$X$.

\begin{defn} \label{quotient} The
quotient variety $X/\Gamma$ is obtained by collapsing each orbit
of $X$ to a point.\end{defn}

If $J$ is a finite group, $c(J)$ denotes the set of conjugacy
classes of $J$.   If $x \in X$,  $\Gamma_x$ denotes the isotropy
group of $x$:
\[
\Gamma_x = \{\gamma \in \Gamma \,:\, \gamma x = x\}.\]

The \emph{extended quotient}, denoted $X\q\Gamma$ is obtained from
$X$ by replacing each orbit $\{\gamma x \,:\, \gamma \in \Gamma\}$
by $c(\Gamma_x)$.  To construct $X\q\Gamma$, we proceed as
follows. Let
\[
\widetilde{X}: = \{(\gamma,x) \in \Gamma \times X\,:\, \gamma x =
x\}.\]

The group $\Gamma$ acts on $\widetilde{X}$ by:
\[
\alpha(\gamma,x) = (\alpha \gamma \alpha^{-1}, \alpha
x)\quad\text{ with $(\gamma, x) \in \widetilde{X}$, $\alpha \in
\Gamma$.}\]

\begin{defn} The \emph{extended quotient}, denoted $X\q \Gamma$, is defined as
\[
X\q\Gamma: = \widetilde{X}/\Gamma\] i.e. $X\q \Gamma$ is the
ordinary quotient (as in definition~\ref{quotient}) for the action
of $\Gamma$ on $\widetilde{X}$.
\end{defn}

The projection
\[
\Gamma \times X \longrightarrow X, \quad (\gamma,x) \mapsto x\]
gives a map
\[\pi\colon
X\q\Gamma \longrightarrow X/\Gamma.\] If $p$ is an orbit in $X$,
i.e. $p = \{\gamma x \,:\, \gamma \in \Gamma \}$, then the
pre-image in $X\q\Gamma$ of $p$ is $c(\Gamma_x)$.  Thus, in
forming $X\q\Gamma$, each orbit $\{\gamma x \,:\, \gamma \in
\Gamma \}$ has been replaced by $c(\Gamma_x)$.

\begin{defn} The map $\pi\colon X\q\Gamma \longrightarrow X/\Gamma$ is \emph{the projection of the extended quotient on the
ordinary quotient.}
\end{defn}
\begin{lem} The projection $\pi\colon X\q \Gamma \to X/\Gamma$ is
a finite morphism of algebraic varieties.
\end{lem}

%The projection
%\[\Gamma \times X \longrightarrow \Gamma, \quad (\gamma,x) \mapsto
%\gamma\] gives a map
%\[X\q\Gamma \longrightarrow c(\Gamma)\]
%This map is not always surjective. If $o(\gamma) = \{\alpha \gamma
%\alpha^{-1}\,:\, \alpha \in \Gamma \}$ is a $\Gamma$-conjugacy class,
%then the pre-image of $o(\gamma)$ in $X\q \Gamma$ is
%$X^{\gamma}/\Cent(\gamma)$ where
Let $X^{\gamma}: = \{x \in X \,:\, \gamma x = x\}$ and denote by
$\Cent(\gamma)$ the $\Gamma$-centralizer of $\gamma$.  We have
\[
X\q\Gamma = \bigsqcup_{\gamma} X^{\gamma}/\Cent(\gamma)\] where
one $\gamma$ is chosen in each $\Gamma$-conjugacy class. Let $e
\in \Gamma$ denote the identity element.   Since $X^e/\Cent(e)$ is
the ordinary quotient $X/\Gamma$, we have $X/\Gamma \subset X\q
\Gamma$.

\begin{defn}  The \emph{reduced quotient} is defined as
\[
(X/\Gamma)_{\rho}: = \pi(X \q \Gamma - X/\Gamma).\]
\end{defn}

\begin{lem}  The reduced quotient is an algebraic variety.
\end{lem}

\begin{lem} The extended quotient is multiplicative: if $\Gamma_1
\times X_1 \to X_1$ and $\Gamma_2 \times X_2 \to X_2$ are as
above, then we have
\[
(X_1 \times X_2)\q(\Gamma_1 \times \Gamma_2) = (X_1\q\Gamma_1)
\times (X_2\q\Gamma_2).\]
\end{lem}

\section{Application to the representation theory of $p$-adic groups}
\label{application} Let $F$ be a local nonarchimedean field, let
$G$ be the group of $F$-rational points in a connected reductive
algebraic group defined over $F$, and let $\Irr(G)$ be the set of
irreducible smooth representations of $G$. We recall the data in
the Bernstein programme: $ \fs\in\fB(G)$, $\fs = [M,\sigma]_G$ is
an inertial class in $G$ (with $M$ a Levi subgroup of $G$ and
$\sigma$ a supercuspidal representation of $M$), $D^{\fs}$ is the
$\Psi(M)$-orbit of $\sigma$ in  $\Irr(M)$, with $\Psi(M)$ the
group of unramified characters of $M$, $W^{\fs} = \{w \in
\Nor_G(M)/M\,:\, w \cdot \fs = \fs\}$ and $D^{\fs}/W^{\fs}$ is the
quotient variety, a component of the Bernstein variety.

We will fix a point $\fs \in \fB(G)$ and write $D = D^{\fs}, W =
W^{\fs}$. Let $\Irr(G)^{\fs}$ denote the $\fs$-component of
$\Irr(G)$ in the Bernstein decomposition of $\Irr(G)$. We will
equip the quotient variety $D/W$ with the Zariski topology, and
$\Irr(G)^{\fs}$ with the Jacobson topology coming from $\Prim \,
\mathcal{H}(G)^{\fs}$. We note that irreducibility is an
\emph{open} condition, and so $(D/W)_{\red}$, the set of reducible
points in $D/W$, is a sub-variety of $D/W$. Let $q$ denote the
cardinality of the residue field of $F$.  Let $\mathcal{H}(G)$
denote the Hecke algebra of $G$, let $\mathcal{H}(G)^{\fs}$ denote
the ideal of $\mathcal{H}(G)$ corresponding to $\fs$ in the
Bernstein decomposition of $\mathcal{H}(G)$, and let $\HP_*$
denote periodic cyclic homology.
\smallskip

\begin{conj} There is an isomorphism
\[\HP_*( \mathcal{H}(G)^{\fs})\cong \rm{H}^*(D\q W)\]
and a continuous bijection
\[
\mu\colon D\q W \to \Irr(G)^{\fs}\] such that:

(1) There is an algebraic family $\pi_t\colon D\q W \to D/W$ of
finite morphisms of algebraic varieties, with $t \in
\mathbb{C}^{\times}$, such that \[\pi_1 = \pi, \quad \quad
\pi_{\sqrt{q}} = (inf.ch.) \circ \mu.\]
%For all $t \in \mathbb{C}^{\times}$, $\pi_t$ restricted to D/W is the identity map of D/W onto itself.
If we let $\mathfrak{X}_t$ be the image of $\pi_t$ restricted to
$D\q W - D/W$ then $\mathfrak{X}_t$ is a flat family of algebraic
varieties such that
\[\mathfrak{X}_1 = (D/W)_{\rho}, \quad \quad
\mathfrak{X}_{\sqrt{q}} = (D/W)_{red}.\]
 (2)  For each irreducible
component $\bc \subset D\q W$ there is a cocharacter
$h_{\bc}\colon \mathbb{C}^{\times} \to D$ such that
\[ \pi_t(x) = \pi(h_{\bc}(t)\cdot x)\] for all $x \in \bc$.  If
$\bc = D/W$ then $h_{\bc} = 1$.
\end{conj}
\begin{thm}  The conjecture is true for $G = \SL(2)$. If $\fs = [T,1]_G$
then $\mathfrak{X}_t$ is the $0$-dimensional variety given by the
polynomial $(x+1)(x-t^2) = 0$.
%D\q W =
%D/W \sqcup\{-1\} \sqcup \{1\}$ with $D = \mathbb{C}^{\times}, W =
%\mathbb{Z}/2\mathbb{Z}$. The generator of $\mathbb{Z}/2\mathbb{Z}$
%sends $z$ to $z^{-1}$.  If $\bc = \{-1\}$ then $h_{\bc} = 1$, if
%$\bc = \{1\}$ then $h_{\bc} = t^2$.
\end{thm}

\section{The general linear group}

\begin{thm}The conjecture is true for $\GL(n)$. \end{thm}
\begin{proof}
The proof uses Langlands parameters, together with some careful
combinatorics. In effect, the $L$-parameters encode the extended
quotient for $\GL(n)$. Let $G = \GL(n) = \GL(n,F)$, $n = mr$,
$\tau$ be an irreducible supercuspidal representation of
$\GL(m,F)$,
\[ \fs = [M,\sigma]_G = [\GL(m)^r, \tau^{\otimes r}]_G.\]
We have \[D = D^{\fs} = (\Cset^{\times})^r, \quad W = W^{\fs} =
S_r.\] Let $W_F$ be the Weil group of $F$, and let $\mathcal{L}_F
= W_F \times \SU(2)$.  Let
$\Hom_{\ssimple}(\mathcal{L}_F,\GL(n,\Cset))$ denote the set of
equivalence classes of Frobenius-semisimple smooth homomorphisms
from $\mathcal{L}_F$ to $\GL(n,\Cset)$.

 For each $n \geq 1$ we have the local
Langlands correspondence \cite{HT}
\[ \rec_F\colon
\Irr(\GL(n,F)) \to \Hom_{\ssimple}(\mathcal{L}_F,\GL(n,\Cset)).\]
We shall write
\[ \Phi(G) =
\Hom_{\ssimple}(\mathcal{L}_F,\GL(n,\Cset)).\]

We will denote by $p$ the following partition of $r$:
\[
a_1 + \cdots + a_1 + \cdots + a_l + \cdots + a_l = r_1a_1 + \cdots
r_la_l = r\] where $a_j$ is repeated $r_j$ times. Let $\gamma \in
S_r$ be the corresponding product of $r_1 + \cdots + r_l$ cycles.
The fixed set $D^{\gamma}$ is a complex torus of dimension $r_1 +
\cdots + r_l$.
%We have
%\[
%D^{\gamma} =
%\left\{(a,\ldots,a,\ldots,b,\ldots,b,\ldots,c,\ldots,c,\ldots,d\ldots,d)\,:\,
%a,b,\ldots,c,d \in \Cset^{\times}\right\}\] where $a,\ldots, b$ is
%repeated $a_1$ times, $c,\ldots, d$ is repeated $a_l$ times. The
%centralizer $\Cent(\gamma)$ is a product of wreath products:
%\[\Cent(\gamma) = \Zset/a_1\Zset \wr S_{r_1} \times \cdots\times
%\Zset/a_l\Zset \wr S_{r_l}.\] The cyclic groups act trivially on
%the fixed set $D^{\gamma}$ and so we have
%\[D^{\gamma}/\Cent(\gamma) = D^{\gamma}/(S_{r_1} \times \cdots \times
%S_{r_l}) = \Sym^{r_1}\Cset^{\times} \times \cdots \times
%\Sym^{r_l}\Cset^{\times}.\]

We recall that $\tau$ is supercuspidal representation of $\GL(m)$.
Now let $\rec_F(\tau) = \eta \in \Irr_m(W_F)$. Denote by $R(j)$
the $j$-dimensional irreducible complex representation of
$\SU(2)$. Corresponding to the partition $p$ we have the
$L$-parameter
\[
\phi = \eta \otimes R(a_1) \oplus \cdots \oplus \eta \otimes
R(a_1) \oplus \cdots \oplus \eta \otimes R(a_l) \oplus \cdots
\oplus \eta \otimes R(a_l)\] where $\eta\otimes R(a_1)$ is
repeated $r_1$ times, $\ldots$, $\eta\otimes R(a_l)$ is repeated
$r_l$ times.

Let $\Psi(W_F)$ denote the group of unramified quasicharacters of
the Weil group $W_F$, and consider the complex torus
$\Psi(W_F)^{r_1+\cdots + r_l}$. The \emph{orbit} of $\phi$ in
$\Phi(G)$, via the action of this complex torus, is
\[
\mathcal{O}(\phi) = \Sym^{r_1}\Cset^{\times} \times \cdots \times
\Sym^{r_l}\Cset^{\times} \subset \Phi(G).\]  Let $\psi_j \in
\Psi(W_F)$ with $1 \leq j \leq r_1 + \cdots + r_l$.   We will map
each $L$-parameter in the orbit $\mathcal{O}(\phi)$ as follows:
\[\psi_1\otimes \eta \otimes R(a_1) \oplus \cdots \oplus \psi_{r_1+
\cdots + r_l}\otimes \eta \otimes R(a_l) \mapsto (\psi_1(\varpi),
\ldots,\psi_{r_1+\cdots + r_l}(\varpi))\in D^{\gamma}\]where
$\varpi$ is a uniformizer in $F$.  This induces a \emph{bijection}
\[
\mathcal{O}(\phi) \cong D^{\gamma}/\Cent(\gamma).\] Let
$\Phi(G)^{\fs}$ denote the $\fs$-component of $\Phi(G)$ in the
Bernstein decomposition of $\Phi(G)$, so that $$\Phi(G)^{\fs} =
\rec_F(\Irr(G)^{\fs}).$$ We now take the disjoint union of the
permutations $\gamma$, one chosen in each $S_r$-conjugacy class.
This creates a \emph{canonical} bijection
\[
\Phi(G)^{\fs} \cong D\q W.\] The reduced quotient $(D/W)_\rho$ is
the hypersurface $\mathfrak{X}_1$ given by the single equation
$\prod_{i \neq j}(z_i - z_j) = 0$. The variety $(D/W)_{red}$ is
the variety $\mathfrak{X}_q$ given by the single equation
$\prod_{i \neq j}(z_i - qz_j) = 0$, according to a classical
theorem \cite[Theorem 4.2]{BZ1}, \cite{Ze}.  The polynomial
equation $\prod_{i\neq j}(z_i - tz_j) = 0$ determines a flat
family $\mathfrak{X}_t$ of hypersurfaces. The hypersurface
$\mathfrak{X}_1$ is the \emph{flat limit} of the family
$\mathfrak{X}_t$ as $t \to 1$, as in \cite[p.77]{EH}.

Let $\bc$ denote the irreducible variety $D^{\gamma}/Z(\gamma)$.
The cocharacter $h_{\bc}$ is given by
\[
h_{\bc}\colon t \mapsto (t^{a_1-1},\ldots,t^{1-a_1},
\ldots,t^{a_r-1},\ldots,t^{1-a_r}) \in D.\]

Finally, we have to use the multiplicativity of the extended
quotient.
\end{proof}

\section{The exceptional group $G_2$}
We have chosen the exceptional group $G_2$ as an awkward example,
requiring many delicate calculations, see \cite{ABP2}.  Let
$\fs=[T,\chi\otimes\chi]_G$ where $T\simeq F^\times\times
F^\times$ is a maximal $F$-split torus of $G=G_2$ and $\chi$ is a
ramified quadratic character of $F^\times$. Let $\{\alpha,\beta\}$
be a basis of a set of roots of $G$ with $\alpha$ short and
$\beta$ long. The group $W^\fs\simeq
\Zset/2\Zset\times\Zset/2\Zset$ is generated by the fundamental
reflections $s_\alpha$ and $s_{3\alpha+2\beta}$, and we have
$D=D^\fs=\left\{\psi_1\chi\otimes\psi_2\chi\,:\,\psi_1,\psi_2\in
\Psi(F^\times)\right\}$. We obtain
%\[
%D^\gamma=\begin{cases}
%\left\{(t,t^{-1}): t \in \Cset^\times\right\},&\text{if $\gamma=
%s_{3\alpha+2\beta}$,}\cr
%\left\{(t,t) : t \in \Cset^\times\right\},&\text{if $\gamma=s_\alpha$,}\cr
%\left\{(1,1), (1,-1), (-1,1), (-1,-1)\right\},&\text{if $\gamma=s_\alpha
%s_{3\alpha+2\beta}$,}
%\end{cases}\]
\[D^\gamma/\Cent(\gamma)=\begin{cases}
\left\{\{(\lambda,\lambda),(\lambda^{-1},\lambda^{-1})\right\}\,:\,
\lambda \in \Cset^\times\},&\text{if $\gamma=s_\alpha$,}\cr
\left\{\{(\lambda,\lambda^{-1}),(\lambda^{-1},\lambda)\}\,:\,
\lambda \in \Cset^\times\right\},&\text{if
$\gamma=s_{3\alpha+2\beta}$,}\cr (1,1) \sqcup (-1,-1) \sqcup
(1,-1),&\text{if $\gamma=s_\alpha s_{3\alpha+2\beta}$.}
\end{cases}\]
Denote by $\mathfrak{C}_1$ the line $x-y=0$ and by
$\mathfrak{C}_2$ the hyperbola $xy - 1=0$.
%A point in the algebraic curve $\mathfrak{C}_1$ has coordinates
%the unordered pair $\{z,z\}$, a point in $\mathfrak{C}_2$ has
%coordinates the unordered pair $\{z^{-1},z\}$ with $z \in
%\Cset^{\times}$.
Setting $\pt_1:=(1,1)$, $\pt_2:=(-1,-1)$, and $\pt_3:=(1,-1)$, we
obtain
\[D\q W = D/W  \sqcup \mathfrak{C}_1 \sqcup \mathfrak{C}_2 \sqcup \pt_1 \sqcup
\pt_2 \sqcup \pt_3, \quad \quad (D/W)_{\rho} = \mathfrak{C}_1 \cup
\mathfrak{C}_2 \cup \pt_3.\]  The cocharacter $h_\bc:
\Cset^{\times} \to D$ is as follows:
\[t\mapsto (t,t^{-1}) \,\text{if $\bc=\mathfrak{C}_1$,}\,
(t^{-1},t^{-1}) \, \text{if $\bc=\mathfrak{C}_2$,}\, (1,t^{-2}) \,
\text{if $\bc=\pt_1$ or $\pt_2$,}\, (1,1)\, \text{if $\bc=\pt_3$.}
\] The variety $\mathfrak{X}_t$ is the union of the line
$x-t^2y=0$, the hyperbola $xy-t^{-2}=0$ and the point $\pt_3$.
This is a flat family.  The $2$ curves admit $2$ intersection
points: $(1,t^{-2})$ and $(-1,-t^{-2})$.  Now let $t = \sqrt q$.
At each of the intersection points,
%given by the polynomial $(x-t^2y)(xy -
%t^{-2}) = 0$ and the isolated point $pt_3$. The $2$ curves in this
%variety admit the $2$ intersection points $\{1,q^{-1}\}$,
%$\{-1,-q^{-1}\}$:
the corresponding parabolically induced representation admits $4$
irreducible inequivalent constituents. As for the point $\pt_3$:
the corresponding parabolically induced representation admits $2$
irreducible inequivalent \emph{tempered} constituents, see
\cite{M}.
\begin{thm} The conjecture is true for the point
$\fs=[T,\chi\otimes\chi]_G$.
\end{thm}
% The Acknowledgements are an un-numbered section
\section*{Acknowledgements} The second author was partially supported by an NSF grant.
% Acknowledgements text here

% etc, etc

% The Appendices part is started with the command \appendix;
% appendix sections are then done as normal sections
% \appendix

% \section{}
% \label{}

% The Acknowledgements are an un-numbered section
%\section*{Acknowledgements} Paul Baum.....
% Acknowledgements text here

%\begin{thebibliography}{00}
% please try to use the bibitem system -
% the references should be in alphabetical order of authors' names.
% Articles with a single author first, author will 1 co-author next,
% then author with several co-authors;

% \bibitem{label}
% Text of bibliographic item

%\bibitem{label}

%\end{thebibliography}

Anne-Marie Aubert, Institut de Ma\-th\'e\-ma\-ti\-ques de Jussieu,
U.M.R. 7586 du
C.N.R.S., Paris, France. Email: aubert@math.jussieu.fr\\

Paul Baum, Pennsylvania State University, Mathematics Department,
University
Park, PA 16802, USA. Email: baum@math.psu.edu\\

Roger Plymen, School of Mathematics, Manchester University,
Manchester M13 9PL,
England. Email: plymen@manchester.ac.uk\\

\end{document}